\def\widebar#1{\overline{#1}}

\def\epsilon{e}
\def\Aut{{\rm Aut}}
\def\Bbb#1{{\bf#1}}
\def\bee{{B}}
\def\ch{{\rm Ch}}
\def\chbar{{\widebar\ch}}
\def\Cor{{\rm Cor}}
\def\End{{\rm End}}
\def\ev{{\rm even}}
\def\gqq{{G'(4)}}
\def\gne{{G(n,\epsilon)}}
\def\gn#1{{G(n,#1)}}
\def\Hom{{\rm Hom}}
\def\mapright#1{\smash{\mathop{\longrightarrow}\limits^{#1}}}
\def\mapdown#1{\Big\downarrow\rlap{$\vcenter{\hbox{$\scriptstyle#1$}}$}}
\def\mne{{M(n,\epsilon)}}
\def\od{{\rm odd}}
\def\pne{{P(n,\epsilon)}}
\def\pra{\par}
\def\proof{{\par\noindent{\bf Proof.}\quad}}
\def\qed{{\quad$\bullet$\par}}
\def\remark{{\par\noindent{\bf Remark.}\quad}}
\def\remarks{{\par\noindent{\bf Remarks.}\quad}}
\def\Res{{\rm Res}}
\def\tee{{\bf T}}
\def\tilg{{\widetilde G}}
\def\tg{{\tilg}}
\def\th{{\widetilde H}}
\def\tm{{\widetilde M}}
\def\tp{{\widetilde P}}
\def\tn{{\widetilde N}}
\def\tx{{\widetilde X}}

\def\emm{M}
\def\enn{N}
\def\pee{P}
\def\eee{E}

\def\alph{{\alpha}}
\def\bet{{\beta}}
\def\delt{{\delta_1}}
\def\deltwo{{\delta_2}}
\def\delthr{{\delta_3}}

\def\ati{1}
\def\bo{2}
\def\burn{3}
\def\lc{4}
\def\eve{5}
\def\huea{6}
\def\hueb{7}
\def\lar{8}
\def\ijlone{9}
\def\ijltwo{10}
\def\ly{11}
\def\lew{12}
\def\pal{13}
\def\rog{14}
\def\rus{15}
\def\san{16}
\def\thomas{17}
\def\wei{18}
\def\yag{19}

\centerline{\bf 3-groups are not determined by their 
integral cohomology rings}
\medskip
\centerline{Ian Leary}
\centerline{Departement Mathematik,}
\centerline{ETH Zentrum,}
\centerline{8092 Z\"urich.}
\medskip

\beginsection Abstract.  

There is exactly one compact 1-dimensional Lie group having 27 components
and nilpotence class three.  We give a presentation for the integral 
cohomology ring of (the classifying space of) this group.  We show that 
the groups of order $3^4$ can be distinguished by their first few 
integral cohomology groups, and exhibit a pair of groups of order $3^5$
having isomorphic integral cohomology rings.  

\beginsection Introduction

The aim of this paper is the study of the integral cohomology rings of
a family of 3-groups.  For each $n\geq 4$ and $\epsilon = \pm 1$ a
member $\gne$ of this family is defined.  The group $\gne$ has order 
$3^n$ and may be presented as follows.  
$$\gne = \langle A,B,C | A^3=B^{3^{n-2}}=C^3=[B,C]=1, [B,A]=C, 
[C,A]= B^{\epsilon 3^{n-3}}\rangle$$
The groups $\gn 1$ and $\gn {-1}$ are not isomorphic to each other.
The main result of this paper (Corollary 14) is that for $n\geq 5$ the integral
cohomology rings of $\gn 1$ and $\gn {-1}$ are isomorphic.  These seem
to be the first examples of $p$-groups for any $p$ having this
property.  An elegant argument due to Alperin and Atiyah 
([\eve], page~86) establishes the existence of groups whose orders
divide by more than one prime having isomorphic integral cohomology rings,
and metacyclic examples are known [\lar].  For $p$-groups there seems 
to be no way to exhibit such groups without actually determining the 
relevant cohomology rings and showing that they are isomorphic, and 
this is what we do.  (See however 
[\ijltwo] for a non-computational method to exhibit $p$-groups having 
isomorphic integral cohomology groups.)  

The groups $\gne$ all occur as normal subgroups of a single compact
1-dimensional Lie group $\tilg$ such that the quotient group is
connected, and our method involves first finding the cohomomology of
the Lie group $\tilg$.  For each prime $p\geq 5$ and $n\geq 4$ two
isomorphism types of groups of order $p^n$ having similar
presentations to those given above for $\gne$ may be defined.  Using
the methods of this paper we have obtained some partial results
concerning the cohomology of these groups.  These partial results have
also been obtained by N.~Yagita using different methods, so here we
merely indicate how they could be proved.  

The first section of the paper consists of a brief discussion of the
method we employ, the second is an examination of the groups $\gne$,
$\tilg$ and various subgroups and quotients, and the third section
contains the main results and proofs.  

\eject

\beginsection Method.

Let $\tee$ be the group of complex numbers of modulus one.  
If $G$ is a finite group and $Z$ a central cyclic subgroup of $G$ then
we define $\tg$ to be the central product of $T$ and $G$ amalgamating
$Z$ with the isomorphic subgroup of $T$ (via some fixed embedding).
The group $\tg$ is a compact Lie group with identity component
isomorphic to $\tee$ and group of components $G/Z$.  The group $G$ is
a normal subgroup of $\tg$ with quotient $T/Z\cong \tee$, and since
$T$ is central in $\tg$, $\bee G$ is a principal $\tee$-bundle over
$\bee\tg$.  Under the natural isomorphism
$H^2(\tg)\cong\Hom(\tg,\tee)$, the first Chern class of this bundle
corresponds to a morphism with kernel $G$.  This construction is
explained in more detail in [\ijlone], and was suggested by
P.~H.~Kropholler and J.~Huebschmann [\huea,\hueb].  

For any group $K$, let $\ch(K)$ be the subring of $H^*(K)$ generated
by Chern classes of complex representations of $K$, and let
$\chbar(K)$ be its ``Mackey closure'', that is the subring of $H^*(K)$
generated by $\ch(K)$ and the images under the transfer of $\ch(H)$ as
$H$ ranges over the finite index subgroups of $K$.  One may ask under
what circumstances either $\ch(K)$ or $\chbar(K)$ is the whole of the
even degree cohomology of $K$.  If $G$ and $\tg$ are as above, then
the following lemma links these properties for $G$ and $\tg$.  

\proclaim Lemma 1.  Let $G$ be a finite group with central cyclic
subgroup $Z$ and construct $\tg$ as above.  Then $\ch(G)=H^\ev(G)$ 
(resp. $\chbar(G)= H^\ev(G)$) if and only if $\ch(\tg)=H^\ev(G)$ (resp.
$\chbar(\tg)=H^\ev(\tg)$) and multiplication by the Chern class of
$\bee G$ as a bundle over $\bee\tg$ is injective on the odd degree
cohomology of $\tg$.  

\proof First we claim that any complex representation of $G$ extends
to one of $\tg$.  This follows from the fact that $Z$ must act via scalar
multiplication in any irreducible $G$-representation, because an
eigenspace for any central element of $G$ is a $G$-summand.  Thus 
the  action may be extended from $Z$ to $T$ so
that $T$ also acts by scalar multiplication, and because the image of
$T$ is central in $\End(V)$ the action extends to $\tg$.  

It follows that the image of $\ch(\tg)$ in $H^*(G)$ is exactly
$\ch(G)$.  Subgroups of $G$ are in one-to-one correspondence with
finite index subgroups of $\tg$.  If $\th$ is a finite index 
subgroup of $\tg$, with $H=\th\cap G$ the corresponding subgroup of 
$G$, then $\Res^\tg_G\Cor^\tg_\th=\Cor^G_H\Res^\th_H$ because 
$G\th=\tg$, and it follows that the image of
$\chbar(\tg)$ is $\chbar(G)$.  From now on the proofs of the
statements concerning $\ch(G)$ and $\chbar(G)$ are identical, so we
consider only the former.  

Now consider the spectral sequence for $\bee G$ as a $\tee$-bundle
over $\bee\tg$.  In this spectral sequence the group $E_3^{n,1}$ is
isomorphic to the cokernel of the map from $H^{n+1}(\tg)$ to
$H^{n+1}(G)$, and is equal to the kernel of multiplication by the
Chern class of the bundle as a map from $H^n(\tg)$ to $H^{n+2}(\tg)$,
so that $H^\ev(\tg)$ maps onto $H^\ev(G)$ if and only if
multiplication by this element is injective on $H^\od(\tg)$.  Now if
$\ch(\tg)=H^\ev(\tg)$, then $\ch(\tg)$ and hence also $\ch(G)$ map
onto $E_3^{\ev,0}$, and so in this case $\ch(G)=H^\ev(G)$ if and only
if $E_3^{\od,1}$ is trivial.  It remains to consider the case when
$\ch(\tg)$ is not the whole of $H^\ev(\tg)$.  In this case pick $x$ of
minimal degree in $H^\ev(\tg)\setminus \ch(\tg)$.  
Since $H^2(\tg)\cong\Hom(\tg,\tee)$ is
contained in $\ch(\tg)$, $x$ has degree at least four.  If it were the
case that $\ch(G)=H^\ev(G)$, then $x$ would have to be congruent to an
element of $\ch(\tg)$ modulo the image of the differential, that is 
$x = z+cy$, where $c$ is the Chern class of $\bee G$ as a bundle over 
$\bee\tg$, and $z$ is an element of $\ch(\tg)$.  But now $y$ and $c$ have
lower degree than $x$, so are also in $\ch(\tg)$, and so we obtain a
contradiction.  \qed

\beginsection The groups $\gne$ and $\tilg$.  

The group $G=\gne$ as presented in the introduction is generated by two
elements, because the element $C$ is already in the subgroup generated
by $A$ and $B$.  It follows that the quotient of $G$ by the
intersection of its maximal subgroups is elementary abelian of rank
two, and hence that $G$ has exactly four maximal subgroups.  
The intersection of these is the subgroup generated by $B^3$ and $C$,
and is isomorphic to $C_{3^{n-3}}\oplus C_3$.  We shall call this
subgroup $\enn$, or $\enn(n,\epsilon)$.  Any element not in $\enn$ is
contained in a unique maximal subgroup.  The maximal subgroup
containing $B$ is isomorphic to $C_{3^{n-2}}\oplus C_3$, and will be
referred to as $\emm$ or $\emm(n,\epsilon)$.  The maximal subgroup
containing $A$ is a non-abelian non-metacyclic group expressible as a
central extension with kernel $C_{3^{n-3}}$ and quotient $C_3\oplus
C_3$ (there is a unique isomorphism type of group having these
properties).  We shall call this subgroup $\pee$.  Note that the
intersection of $\pee$ and $\emm$ is $\enn$.  The other two maximal
subgroups, those containing $AB$ and $AB^2$, are non-abelian
metacyclic groups with a cyclic subgroup of index three, except in the
case of $G(4,-1)$ when they are isomorphic to $\pee(4,-1)$, which is
the non-abelian group of order 27 and exponent 3. 

The centre of
$\gne$ is cyclic of order $3^{n-3}$ generated by $B^3$.  The quotient
of $\gne$ by its centre shall be called $\eee(n,\epsilon)$, and is
the (unique) non-abelian group of order 27 and exponent 3.  Since any
order~3  normal subgroup of a 3-group must be central, it follows
that $\gne$ has a unique such subgroup, the subgroup generated by
$B^{n-3}$.  The quotient groups $G(n,1)/\langle B^{n-3}\rangle$ and 
 $G(n,-1)/\langle B^{n-3}\rangle$ are isomorphic, via a map sending
the elements $A$, $B$ and $C$ in one group to the elements of the same
name in the other, and it follows that the lattices of normal
subgroups of $G(n,1)$ and $G(n,-1)$ are isomorphic.  If $n$ is at
least~5 then a far stronger statement is true.  In this case, if we write
elements of $\gne$ as words of the form $A^iB^jC^k$ then a
collection of elements forms a normal subgroup of $G(n,1)$ if and 
only if the 
`same' collection of elements forms a normal subgroup of $G(n,-1)$,
and these two subgroups are isomorphic provided that they are proper.
The author can think of no proof of this fact apart from direct
calculation in the maximal subgroups.  

It is known that the isomorphism type of a $p$-group is determined by
that of its integral group ring [\rog,\wei], and R.~Sandling has shown
the author a proof that the mod-3 group rings of the groups $G(n,1)$
and $G(n,-1)$ are not isomorphic [\san].  
The group $\gne$ has only 1- and 3-dimensional irreducible
representations because it has an abelian subgroup (the subgroup
$M$) of index three.  The sizes of conjugacy classes in $G(n,1)$ and
$G(n,-1)$ are identical, but the character tables of $G(n,1)$ and
$G(n,-1)$ are different.  The character table for $\gne$ contains the
entry $\eta(2+\eta^{\epsilon 3^{n-3}})$ for each primitive $3^{n-2}$th
root of unity $\eta$, while the character table for $G(n,-\epsilon)$
does not.  This gives a proof that $G(n,1)$ and $G(n,-1)$ are not
isomorphic.  (It is also reasonably straightforward to prove this fact
directly.)  The following proposition describes the $\Lambda$-ring
structure of the representation ring of $\gne$.  

\def\thrr{{3^{n-4}}}

\proclaim Proposition 2.  The representation ring of $\gne$ is
generated by $\theta$, $\psi$ of dimension one and $\chi$, $\bar\chi$,
$\xi$, $\bar\xi$ of dimension three (where the bar indicates the
dual representation) subject to the following relations, together
with commutativity and the relations implied by $\widebar{XY}=\bar X
\bar Y$.  
$$\theta^3=1,\quad \psi^{3^{n-3}}=1,\quad \theta\chi=\chi,\quad
\psi^\thrr\chi=\chi,$$
$$\chi^2=3\bar\chi,\quad\chi\bar\chi
=(1+\theta+\theta^2)(1+\psi^\thrr+\psi^{-\thrr}),$$
$$\xi\chi=\xi\bar\chi=\xi(1+\psi^\thrr+\psi^{-\thrr}),$$
$$\theta\xi=\xi,\quad\xi\bar\xi=\chi+\bar\chi+1+\theta+\theta^2,$$
$$\xi^2=\bar\xi\psi(1+2\psi^{\epsilon\thrr}).$$
The $\Lambda$-ring structure is given by the following equations. 
$$\Lambda^2(\chi)=\bar\chi\psi^3,\quad\Lambda^3(\chi)=\psi^3,$$
$$\Lambda^2(\xi)=\bar\xi\psi^{(1+\epsilon\thrr)},\quad
\Lambda^3(\xi)=\psi^{(1+\epsilon\thrr)}.$$
\par
\proof Direct calculation.  In the above statement, $\theta$, $\psi$
and $\chi$ come from representations of the quotient of $\gne$ by its
central subgroup of order three, and $\xi$ is a faithful
representation of $\gne$ obtained by inducing up a representation of
$\emm(n,\epsilon)$ that is faithful on $\langle B\rangle$.  \qed

I do not know if the representation rings of $G(n,1)$ and $G(n,-1)$
are isomorphic.  I conjecture that even after quotienting out by the
ideal generated by 9 the rings are not isomorphic.  This would imply 
by Atiyah's theorem [\ati] that the $K$-theory ring can distinguish
the groups $G(n,1)$ and $G(n,-1)$.  For $n$ at least 5, the
representation rings modulo the ideal generated by 3 are isomorphic,
via the map that sends $\xi$ to $-\xi\psi^{2.3^{n-5}}$, $\bar\xi$ to 
$-\bar\xi\psi^{7.3^{n-5}}$, and fixes the other generators.  This
isomorphism commutes with taking duals but not with the
$\Lambda$-ring structure.  

One might hope to distinguish the representation rings of $G(n,1)$ 
and $G(n,-1)$ by using their $\Lambda$-ring structure.  For example, 
Grothendieck has defined a filtration (the $\gamma$-filtration) on any 
augmented $\Lambda$-ring, and one may consider the representation rings
of the groups modulo the layers of this filtration.  A theorem of 
Atiyah [\ati] and the fact that Chern classes generate 
the even degree cohomology
of $\gne$ (see Corollary~9) imply that this is equivalent to studying 
the image of $K^0(\bee\gne)$ in $K^0$ of its skeleta.  The layers of 
the $\gamma$-filtration for $\gne$ may be computed directly from
Proposition~3, but the task is simplified by also using information
concerning the low degree cohomology of $\gne$ as provided by
Theorem~13.  The author has been able to show that the representation
rings of $G(5,1)$ and $G(5,-1)$ are isomorphic after quotienting out
the third layer of the $\gamma$-filtration, or equivalently that the
images of $K^0(\bee G(5,1))$ and $K^0(\bee G(5,-1))$ in $K^0$ of their
respective 5-skeleta are isomorphic.  

There is one other finite group closely related to the groups $\gne$.
We may view $\gne$ as an extension with kernel $\emm(n,\epsilon)$ and
quotient cyclic of order three.  In each case $\gne$ is the
corresponding split extension.  It is easily verified that the second
cohomology group of $C_3$ with coefficients in the module
$\emm(n,\epsilon)$ is trivial, unless $n = 4$ and $\epsilon = -1$, in
which case it has order three.  This gives rise to another group of
order $3^4$, which we shall call $\gqq$.  This group may be presented
as follows.  
$$\gqq= \langle A,B,C | B^9=C^3=[B,C]=1, [B,A]=C,
[C,A]=B^{-3}=A^3\rangle$$
The character table (and hence also the representation ring) of this
group is identical to that of $G(4,-1)$, although the $\Lambda$-ring
structure is slightly different.  

Each of the groups $\gne$ and $\gqq$ is nilpotent of class three, 
and has cyclic centre of index~27 with quotient group isomorphic to
$E$, the non-abelian group of order 27 and exponent 3.  Now let $G$ be
any of the groups $\gne$ or $\gqq$, and apply the construction of
section~1 to form a Lie group $\tg$ as the central product of $G$
and $\tee$ amalgamating the centre of $G$ with the isomorphic subgroup
of $\tee$.  In each case the resulting group is nilpotent of class
three, has identity component isomorphic to $\tee$ and group of
components $E$.  There is however, only one isomorphism type of group
having these properties, because the action of $\Aut(E)$ on
$H^2(E;\tee)\cong H^3(E;\Bbb Z)$ is transitive on non-zero elements,
and $\tee\times E$ is nilpotent of class two.  Thus the construction
of section~1 allows us to embed all of the groups $\gne$ and $\gqq$
into a single Lie group $\tg$ as normal subgroups with connected
quotient.  This Lie group may be presented as follows, where $\tee$ is
considered to be a subgroup of the complex numbers of modulus one, and
$\omega$ is $\exp(2\pi i/3)$.
$$\tg=\langle X,Y,Z,\tee|X^3=Y^3=Z^3=1, \hbox{$\tee$ central, }
[Y,Z]=1, [Y,X]=Z,[Z,X]=\omega\rangle$$
This presentation enables one to express any element of $\tg$ in the 
form $X^iY^jZ^kt$ for some $t\in \tee$.  
It is easy to check that if $\eta=\exp(2\pi i/3^{n-2})$, then the
subgroup with generators $X$, $Y\eta^\epsilon$ and $Z$ is isomorphic
to $\gne$ with generators $A$, $B$ and $C$, and we fix this embedding
from now on.  The group $\tg$ has four subgroups of index three, one
of which is abelian, and three of which are the unique non-abelian Lie
group consisting of nine circles.  We shall refer to the subgroup
generated by $\tee$, $Y$ and $Z$ as $\tm$, this being the abelian
subgroup of index three, and the subgroup generated by $\tee$, $X$ and
$Z$ shall be called $\tp$.  Regarding the embedding of $\gne$ in $\tg$
as an inclusion the following equalities hold.  
$$\pee(n,\epsilon)=\tp\cap\gne,\quad\emm(n,\epsilon)=\tm\cap\gne$$

In the next section we shall require a description of some elements of
$\Hom(\tg,\tee)$ having kernel isomorphic to $\gne$ or $\gqq$.  It is
also possible to classify the isomorphism types of 3-groups that can
occur as kernels of maps from $\tg$ to $\tee$, and so we combine these
statements in the following proposition.  Before stating the proposition, 
it is convenient to make the following definition, which we shall use 
frequently in the sequel.
  
\proclaim Definition.  For any ring $R$ and any $R$-module $M$, we say 
that a subset $S$ is a basis for $M$ if zero is not an element of $S$
and $M$ is isomorphic to the direct sum of the submodules generated by 
the elements of $S$.  Note that any finitely generated module for a 
principal ideal domain has such a basis.  

\proclaim Proposition 3.  Define $\alph$ (resp. $\bet$, $\delt$) in 
$\Hom(\tg,\tee)$ by insisting that it maps $X^iY^jZ^kt$ to $\omega^i$
(resp. $\omega^j$, $t^3$).  Then $\alph$ and $\bet$ have order three, 
$\delt$ has infinite order, and these elements form a basis for the 
group $\Hom(\tg,\tee)$.  
The elements of $\Hom(\tg,\tee)$ having
kernel of order $3^n$ are those of the form $\pm3^{n-4}\delt+\gamma$,
where $\gamma$ is in the subgroup spanned by $\alph$ and $\bet$.  The
action of a general automorphism of $\tg$ on $\Hom(\tg,\tee)$ is as
follows, where $m$ is either 0 or 1, $i$ is either 1 or 2, and $j$ is
0, 1, or 2.  
$$\eqalign{\delt&\mapsto(-1)^m\delt-j\alph\cr
\bet&\mapsto(-1)^m\bet+j\alph\cr
\alph&\mapsto i\alph}$$
It follows that as orbit representatives among elements having kernel
a 3-group we may take $\thrr\delt$, $\thrr\delt+\bet$,
$\thrr\delt-\bet$ and $\delt+\bet+\alph$.  The first of these has
kernel containing a $(C_3)^3$ subgroup, while the others correspond to
$G(n,-1)$, $G(n,1)$ and $\gqq$ respectively.  

\proof The only portion of the statement that we actually need in the
sequel is the last sentence, which can be checked by a simple
calculation, so we omit the proof.  \qed

\remark Instead of using $\gne$ to construct $\tg$, one could start
from $\tg$ and use Proposition~3 to define $\gne$.  The obligation to
provide a proof for Propostion~3 would then be greater, and one would
still have to write down presentations for the groups $\gne$, so this
approach would not save any labour.  

\beginsection Cohomology. 

The cohomology of the Lie group $\tg$ will be calculated below, using
the spectral sequence for $\tg$ as an extension with kernel $\tp$ and
quotient cyclic of order three.  From this the cohomology of $\gne$
can be calculated easily.  To minimise the number of letters employed
to represent cohomology classes we adopt the following convention.

\proclaim Notation.  If $\xi$ represents an element of $H^*(\tg)$, 
we shall use the same symbol to represent its image in $H^*(\tp)$,
$H^*(\tm)$, or $H^*(\tn)$.  In case of ambiguity we shall refer to the
element $\xi$ of $H^*(\tx)$ as $\xi(\tx)$.  If we define elements
$\xi(\tp)$ and $\xi(\tm)$, we do not wish to imply that these elements
are images of an element $\xi(\tg)$, but merely that their images in
$H^*(\tn)$ are equal.  This convention extends in the obvious way to
other subgroups of $\tg$.  

We now define generators for $H^*(\tm)$, and use them to define
elements of $H^*(\tg)$ which will later be shown to form a generating
set.

\proclaim Proposition 4.  Define an element $\tau$ (resp.\ $\beta$,
$\gamma$) in $\Hom(\tm,\tee)\cong H^2(\tm)$ by insisting that it maps 
$Y^jZ^kt$ to $t$ (resp. $\omega^j$, $\omega^k$), and let $\mu$ be any
non-zero element of $H^3(\tm)$.  Then $H^*(\tm)$ is generated by
$\tau$, $\beta$, $\gamma$ and $\mu$, subject only to the relations 
$3\beta = 0$, $3\gamma=0$, $3\mu =0$ (and of course any relations
implied by anticommutativity).  The action of conjugation by $X$ on
$H^*(\tm)$ sends $\tau$ to $\tau+\gamma$, sends $\gamma$ to
$\gamma+\beta$, and fixes $\beta$ and $\mu$.  
The restriction map from $H^*(\tm)$ to
$H^*(\mne)$ is surjective and has kernel the ideal generated by
$3^{n-3}\tau-\epsilon\beta$.  

\proof The cohomology of $\tm$ is easily shown to be as claimed 
(note that $\tm\cong \tee\oplus C_3\oplus C_3$), as is
the action of conjugation by $X$.  In the Gysin sequence for $\bee M$
as a $\tee$-bundle over $\bee\tm$ the map from $H^n(\tm)$ to
$H^{n+2}(\tm)$ is multiplication by the image in $H^2(\tm)$ of
$3^{n-4}\delt-\epsilon\bet$, which is $3^{n-3}\tau-\epsilon\beta$.
The kernel of multiplication by this element is trivial, which implies
that the restriction map from $H^*(\tm)$ to $H^*(M)$ is surjective,
and its kernel is clearly as claimed.  \qed

\proclaim Definition/Proposition 5.  In addition to the elements
$\alph$, $\bet$ and $\delt$ of $H^2(\tg)$ defined in the statement of
Proposition~3, define elements $\deltwo$, $\delthr$, $\zeta$ and $\mu$
of $H^*(\tg)$ as follows.  Let $\rho$ be a 3-dimensional irreducible
representation of $\tg$ whose restriction to $\tm$ contains the
1-dimensional representation with first Chern class $\tau$, and define
$\delthr$ to be the third Chern class of $\rho$.  Define $\deltwo$ and
$\zeta$ as transfers from $\tm$ by the following equations.  
$$\deltwo= \Cor^\tg_\tm(\tau^2)\qquad
\zeta=\Cor^\tg_\tm(\tau^2(\beta+\gamma))=\deltwo\beta(\tg)+
\Cor^\tg_\tm(\tau^2\gamma)$$
The element $\mu(\tg)$ may be uniquely defined by requiring that its
image in $H^3(\tm)$ is $\mu(\tm)$.  The following equation relates 
$\delt(\tg)$ and $\tau(\tm)$.  
$$\delt=\Cor^\tg_\tm(\tau)-\beta$$

\proof The last two sentences of the above statement require a proof,
the rest being definitions.  The assertion concerning $\mu(\tg)$ is
equivalent to the assertion that $H^3(\tg)$ has order three and maps
injectively to $H^3(\tm)$.  This can be shown easily 
by considering the spectral
sequence for $\tg$ as an extension with kernel $\tee$ and quotient
$E$, but will also follow from our study of the spectral sequence for
$\tg$ as an extension with kernel $\tp$ and quotient $C_3$, so we
postpone the proof of this statement until the end of the proof of
Lemma~8.  

To verify the equation relating $\Cor^\tg_\tm(\tau)$ to $\delt$, note
that the transfer from $H^2(\tm)$ to $H^2(\tg)$ is equal to the
following composite, 
$$H^2(\tm)\cong\Hom(\tm,\tee)\mapright{t^*}\Hom(\tg_{{\rm ab}},\tee)
\cong H^2(\tg),$$
where $t^*$ is the map induced by the classical transfer map from
$\tg_{{\rm ab}}$ to $\tm$.  Using this we may describe $\Cor(\tau)$ as
an element of $\Hom(\tg,\tee)$ by the following equations.  
$$\eqalign{\Cor(\tau)(X^iY^jZ^kt)&=\tau(Y^jZ^kt)\tau(XY^jZ^ktX^{-1})
\tau(X^{-1}Y^jZ^ktX)\cr
&=\tau(Y^jZ^kt)\tau(Y^jZ^{j+k}t\omega^k)\tau(Y^jZ^{2j+k}t\omega^{2k+j})\cr
&=t^3\omega^j\cr
&=(\delt+\bet)(X^iY^jZ^kt)\cr}$$
\qed


\proclaim Theorem 6.  Define $\gamma\in \Hom(\tp,\tee)$ by the
equation $\gamma(X^iZ^kt)=\omega^k$, and let $\alph$, $\delt$,
$\deltwo$, $\delthr$ be the restrictions to $\tp$ of the elements of
$H^*(\tg)$ having the same names.  Then these five elements generate
$H^*(\tp)$ subject only to the following relations.  
$$3\gamma=3\alph=0,\qquad \alph^3\gamma=\gamma^3\alph,$$
$$\alph\delt=0=\gamma\delt,\qquad \alph\deltwo=0,\qquad\gamma\deltwo
=-\gamma^3+\alph^2\gamma,$$
$$\delt^2=3\deltwo,\qquad \delt\deltwo=9\delthr,\qquad
\deltwo^2=3\delthr\delt+\gamma^4-\alph^2\gamma^2.$$
The action of conjugation by $Y$ on $H^*(\tp)$ sends $\gamma$ to
$\gamma-\alph$ and fixes the other generators.  
\pra
The image of the map from $H^*(\tp)$ to $H^*(\pne)$ is the whole of 
$H^{{\rm even}}(\pne)$, and the kernel is the ideal generated by
$3^{n-4}\delt$.  As a module for $H^{{\rm even}}(\pne)$, $H^{{\rm
odd}}(\pne)$ is generated by two elements $\mu_1$ and $\mu_2$ of
degree three, subject to the following relations.  
$$3\mu_1=3\mu_2=0,\qquad\mu_1\delt=\mu_2\delt=0,\qquad
\mu_1\gamma=\mu_2\alph,$$
$$\mu_1\deltwo=0,\qquad\mu_2\deltwo=-\gamma^2\mu_2-\alph^2\mu_2,
\qquad \alph^3\mu_2=\gamma^3\mu_1.$$
The ring structure of $H^*(\pne)$ is determined by the above relations
together with the relation 
$$\mu_1\mu_2=\cases{0&for $n>4$\cr 3\delthr&for $n=4$.}$$

\proof The integral cohomology of $\tp$ and $\pne$ is computed in
[\ijlone], where the group $\pne$ is called $P(n-1)$, and a 
different generator is taken in degree four from the one used above. 
Assuming the statements contained in [\ijlone], we only need to
compare the restrictions to $\tp$ of our generators for $H^*(\tg)$ and
the generators taken for $H^*(\tp)$ in [\ijlone].  The representation
used to define $\delthr(\tg)$ restricts to $\tp$ as the representation
used there to define the degree six generator in $H^*(\tp)$.  The
behaviour of the degree two generators defined as homomorphisms when
restricted to $\tp$ is clear.  For the generators defined in terms of
the transfer the required relations follow from the observation that
$\tm\tp=\tg$, and so the double coset formula gives the equation
$\Cor^\tp_\tn\Res^\tm_\tn=\Res^\tg_\tp\Cor^\tg_\tm$.  
\qed

\proclaim Proposition 7.  The restrictions to $\tm$ and $\tp$ of the
elements of $H^*(\tg)$ of Definition~5 are either of the form
``element maps to element having the same name'', or are included in
the following lists.  
\pra
Restrictions to $\tp$:
$$\Res(\bet)=0,\qquad\Res(\mu)=0,\qquad\Res(\zeta)=\alph^2\gamma-\gamma^3.$$
\pra
Restrictions to $\tm$:
$$\Res(\alph)=0,\qquad\Res(\delt)=3\tau,\qquad
\Res(\zeta)=\beta^2\gamma-\gamma^3,$$
$$\Res(\deltwo)=3\tau^2-\tau\beta-\gamma^2+\gamma\beta+\beta^2,\qquad
\Res(\delthr)=\tau^3+\tau^2\beta-\tau\gamma^2+\tau\gamma\beta.$$

\proof For the elements defined using the transfer, the double coset
formula suffices to obtain the above results.  As an example, the
restriction to $\tp$ of $\zeta$ may be found as follows.  
$$\eqalign{\Res^\tg_\tp(\zeta)
&=\Res^\tg_\tp\Cor^\tg_\tm(\tau^2(\beta+\gamma))\cr
&=\Cor^\tp_\tn\Res^\tm_\tn(\tau^2(\beta+\gamma))\cr
&=\deltwo(\tp)\gamma(\tp)=\alph^2\gamma-\gamma^3.}$$
The restriction to $\tm$ of the representation $\rho$ used to define
$\delthr$ contains a summand with Chern class $\tau$, so must also
contain summands with Chern classes $\tau+\gamma$ and
$\tau-\gamma+\beta$, these being images of $\tau$ under the action of
conjugation by powers of $X$ on $H^*(\tm)$.  Its third Chern class
(which is by definition the restriction to $\tm$ of $\delthr$) is the
product of these three elements.  \qed

\proclaim Lemma 8.  The spectral sequence with integer coefficients
for $\tg$ as an extension with kernel $\tp$ and quotient of order
three collapses.  The seven elements of $H^*(\tg)$ of Definition~5
generate the $E_2$-page.  The elements $\alph$, $\delt$, $\deltwo$,
$\delthr$ and $\zeta$ yield elements in $E_2^{0,*}$, $\beta$ yields
a generator for $E_2^{2,0}$, and $\mu$ yields a generator for
$E_2^{1,2}$.  The map of spectral sequences induced by the diagram 
$$\matrix{\tn&\mapright{}&\tm&\mapright{}&C_3\cr
\mapdown{}&&\mapdown{}&&\mapdown{}\cr
\tp&\mapright{}&\tg&\mapright{}&C_3}$$
is injective on the $E_2^{i,j}$ such that $i+j$ is odd.  
The ring structure of the $E_2$-page is given by the
following relations.  
$$3\alph= 3\beta=0,\qquad 3\mu=0,\qquad 3\zeta=0,$$
$$\alph\delt=0,\qquad\delt^2=3\deltwo,\qquad\alph\deltwo=0,
\qquad\delt\deltwo=9\delthr,$$
$$\alph\zeta=0,\qquad\delt\zeta=0,\qquad\alph^2\beta=0,$$
$$\alph\mu=0,\qquad\delt\mu=0,\qquad 27\delthr^2-\deltwo^3=\zeta^2.$$

\proof First we split $H^*(\tp)$ as a sum of indecomposable modules
for the action of $C_3$.  From the relations given in Theorem~6 it is
easy to show that the elements $\delthr^i\delt$, $\delthr^i\deltwo$, 
$\delthr^i\gamma^j$, $\delthr^i\gamma^j\alph$,
$\delthr^i\gamma^j\alph^2$, $\delthr^i\alph^{j+3}$ form a basis for
$H^*(\tp)$, where $i$ and $j$ are any positive integers.  (Recall that
we have defined a basis for an abelian group to be a set of elements
not containing the identity element such that the group is equal to
the direct sum of the cyclic subgroups generated by those elements.)  
The $C_3$-submodules spanned by $\delthr^i$, $\delthr^i\delt$ and
$\delthr^i\deltwo$ are direct summands isomorphic to the trivial
$C_3$-module $\Bbb Z$, and the monomials of fixed degree in $\delthr$
and fixed total degree in $\alph$ and $\gamma$ (there are at most four
such for any choice of degrees) form a $C_3$-summand.  The elements
$\delthr^i\alph$ and $\delthr^i\gamma$ form an indecomposable
$C_3$-summand with underlying group $C_3\oplus C_3$, and the elements 
$\delthr^i\alph^2$, $\delthr^i\alph\gamma$, $\delthr^i\gamma^2$ form
an indecomposable $C_3$-summand which must be a free $\Bbb
F_3C_3$-module of rank one.  For $j\geq 3$ and for any $i$, the 
$\Bbb F_3C_3$-module generated by the four monomials of degree $i$ in
$\delthr$ and degree $j$ in $\alph$ and $\gamma$ splits as a direct
sum of a trivial module generated by
$\delthr^i(\gamma^j-\gamma^{j-2}\alph^2)$ and an indecomposable module
containing $\delthr^i\alph^j$, $\delthr^i\gamma^{i-1}\alph$ and
$\delthr^i\gamma^{j-2}\alph^2$, which is therefore free.  

Using the above $C_3$-splitting of $H^*(\tp)$, it is easy to check
that the elements $\delthr^i\delt$, $\delthr^i\deltwo$,
$\delthr^i\alph^j$ and $\delthr^i(\gamma^{j+3}-\gamma^{j+1}\alph^2)$
where $i,j\geq 0$ form a basis for the fixed point subring.  
Identifying elements of $H^*(\tg)$ with their images in $H^*(\tp)$ and
applying Theorem~6 and Proposition~7, we see that for $i\geq 0$ the
following equalities hold.  
$$\eqalign{\gamma^{3i+3}-\gamma^{3i+1}\alph^2&=\zeta^{i+1}\cr
\gamma^{3i+4}-\gamma^{3i+2}\alph^2&=\zeta^i(\deltwo^2-3\delthr\delt)\quad
(=\zeta^i\deltwo^2\hbox{ if } i>0)\cr
\gamma^{3i+5}-\gamma{3i+3}\alph^2&=-\zeta^{i+1}\deltwo}$$
It is now easy to see that the elements $\delthr^i\delt$,
$\delthr^i\deltwo$, $\delthr^i\alph^j$, $\delthr^i\zeta^{j+1}$,
$\delthr^i\zeta^{j+1}\deltwo$, $\delthr^i(\deltwo^2-3\delthr\delt)$ and
$\delthr^i\zeta^{j+1}\deltwo^2$ where $i,j \geq 0$ also form a basis
for the fixed point subring.  This already implies that all
differentials in the spectral sequence are trivial on $E^{0,*}$.  It
is now easy to see that the relations claimed between the elements
$\alph$, $\delt$, $\deltwo$, $\delthr$ and $\zeta$ are exactly the
relations that do hold between them as elements of $H^*(\tp)^{C_3}$,
or equivalently as elements of $E_2^{0,*}$.  Note that the fact that 
$\deltwo^2-3\delthr\delt$ has order three follows from the 
given relations by expanding
$\delt^2\deltwo$ in two different ways.    

It is clear that in the spectral sequence $\beta$ yields an element of
$E_2^{2,0}$ because as a homomorphism from $\tg$ to $\tee$, $\beta$ has
kernel $\tp$.  Cup product with this element of $E_2^{2,0}$ gives a
surjection from $E_2^{0,j}$ to $E_2^{2,j}$ and an isomorphism from
$E_2^{i+1,j}$ to $E_2^{i+3,j}$ for all $i$.  We now consider the
cohomology of the quotient $C_3$ with coefficients in the various
modules that occur in our decomposition of $H^*(\tp)$.  Of the types
of module that occur, each has second cohomology group of order 3
except for the free $\Bbb F_3C_3$-module, which has trivial positive
degree cohomology.  The subset of the above basis of $H^*(\tp)^{C_3}$
corresponding to free $\Bbb F_3C_3$-summands of $H^*(\tp)$ consists of
the elements $\zeta^i\alph^{j+2}$, so our relations between the
elements of $E_2^{*,*}$ of even total degree may be completed by
adding the relations $3\beta=0$ and $\alph^2\beta=0$.  

Now we consider the elements of $E_2^{*,*}$ of odd total degree.  The
first cohomology group of $C_3$ with coefficients in $\Bbb Z$ or in
$\Bbb F_3C_3$ is trivial, and the corresponding group with
coefficients in each of the other two types of indecomposable summands
of $H^*(\tp)$ is cyclic of order three.  Using the above
$C_3$-splitting of $H^*(\tp)$ it is now easy to find the dimension
over $\Bbb F_3$ of $E_2^{1,j}$.  If we define $P(t)$ to be the power
series whose coefficient of $t^j$ is the dimension of $E_2^{1,j}$,
then the following equation describes $P(t)$.  
$$P(t)=(t^6-t^4+t^2)/(1-t^6)(1-t^2)$$
It is reasonably easy to describe the rest of the multiplicative
structure of the $E_2$-page (that is, the products involving at least
one element of odd total degree) directly, but this will follow from
the assertion concerning restriction to the spectral sequence for
$\tm$ expressed as an extension with kernel $\tn$ and quotient of
order three.  Let $\widebar E_*^{*,*}$ be the spectral sequence for
this extension.  Then $\widebar E_*^{*,*}$ collapses, and the
$E_2$-page is isomorphic as a ring to $H^*(\tm)$, where $\tau$ and
$\gamma$ yield elements of $\widebar E_2^{0,2}$, $\beta$ yields an
element of $\widebar E_2^{2,0}$, and $\mu$ yields an element of
$\widebar E_2^{1,2}$.  It is easy to check that $E_2^{1,2}$ maps
isomorphically to $\widebar E_2^{1,2}$ using the cohomology long exact
sequence associated to the following short exact sequence of
$C_3$-modules, where $K$ (which is defined by this sequence) is a
trivial module isomorphic to the integers modulo three.  
$$0\longrightarrow K\longrightarrow H^2(\tp) \longrightarrow
H^2(\tn) \longrightarrow 0$$
The images of $\alph$, $\delt$, $\deltwo$, $\delthr$ and $\zeta$ in
$\widebar E_2^{0,*}$ are just the coefficients of $\beta^0$ in their
restrictions to $\tm$ (see Proposition~7).  In $\widebar E_2^{*,*}$,
$\alph=0$, and $\delt$ is divisible by three, so here the relations
$\alph\mu=0$ and $\delt\mu=0$ hold.  The subring of $E_2^{0,*}$
generated by $\deltwo$, $\delthr$ and $\zeta$ maps injectively to
$\widebar E_2^{0,*}$ (as can be seen by checking the images of the
basis for this ring given earlier), and we shall temporarily refer to
this subring of $\widebar E_2^{0,*}$ as $R$.  The $R$-module generated
by $\mu$ is isomorphic to $R/3R$, and it may be checked that the
Poincar\'e series for this subgroup of $\widebar E_2^{1,*}$ is equal
to the series $P(t)$ above.  From this it follows that the odd total
degree subgroups of $E_2^{*,*}$ map injectively to $\widebar
E_2^{*,*}$, that no other generators are required for $E_2$, and that
the relations in $E_2$ involving $\mu$ are as claimed.  

To show that the spectral sequence $E_2^{*,*}$ collapses we must show
that each of the seven generators survives.  From Definition~5 it is
clear that each of the even degree generators actually comes from an
element of $H^*(\tg)$, so must survive.  The definition given there
for $\mu$ is more nebulous, and it has not yet been shown that such an
element exists, although the preceeding paragraph shows that such an
element is unique.  The only way for $\mu\in E_2^{1,2}$ to fail to
survive is for $d_3(\mu)$ to be a non-zero multiple of $\beta^2$, and
there are many ways to see that this cannot happen.  For example, in
$\widebar E_*^{*,*}$, $d_3(\mu)$ is zero, but $E_3^{4,0}$ maps
injectively to $\widebar E_3^{4,0}$, so $d_3(\mu)$ must be zero in
$E_*^{*,*}$ too.  Alternatively one may use the fact that for any
split extension and any trivial coefficients, no differential in the
spectral sequence can hit the base-line (see [\lew], where the
condition that the coefficients be trivial is omitted).  This
completes the proofs of Lemma~8 and Definition~5.  \qed

\proclaim Corollary 9.  Chern classes of representations generate the
even degree cohomology of $\gne$ and $\gqq$.  

\proof By choosing a slightly different set of generators for
$H^*(\tg)$ from the one we have chosen, and using the arguments of
Lemma~8, it is easy to show that Chern classes generate the even
degree cohomology of $\tg$.  It is also known that Chern classes
generate the even degree cohomology of the subgroups $M(n,\epsilon)$
and the abelian maximal subgroup of $\gqq$.  Applying one implication
of Lemma~1 to $\tm$, it may be seen that multiplication by the Chern 
class of the $\tee$-bundle $\bee M(n,\epsilon)$ over $\bee \tm$ 
is injective on the odd
degree cohomology of $\tm$.  However, it was shown in Lemma~8 that the
odd degree cohomology of $\tg$ maps injectively to that of $\tm$, and
it follows that multiplication by the Chern class of the $\tee$-bundle
$\bee\gne$ (resp. $\bee\gqq$) over $\bee\tg$ is injective on the odd
degree cohomology of $\tg$.  Now the other implication of Lemma~1
applied to $\tg$ gives the required result.  \qed

\remark In the introduction we remarked that there are two groups of
order $p^n$ for each prime $p\geq 5$ and each $n\geq4$ having similar
presentations to the groups $\gne$.  These groups occur as normal
subgroups of the (unique) Lie group having $p^3$ circular components
and nilpotence class three, just as $\gne$ occurs within $\tg$. This
Lie group is expressible as a split extension with kernel the unique
non-abelian Lie group having $p^2$ circular components, and quotient
of order $p$.  I do not know whether the spectral sequence for this
extension collapses, but the methods of Lemma~8 can be used to show
that for this Lie group corestrictions of Chern classes generate the
even degree cohomology, and that the odd degree cohomology maps
injectively to the cohomology of its (unique) abelian maximal
subgroup.  As in Corollary~9 it is possible to deduce that for the
corresponding finite groups, corestrictions of Chern classes generate
the even degree cohomology.  This result has been obtained using other
methods by N. Yagita [\yag].  It may also be shown that for these groups
Chern classes alone do not suffice to generate the even degree
cohomology [\ly].  

\proclaim Theorem 10. The integral cohomology of the Lie group $\tg$
is generated by the seven elements of Definition~5 subject to the
following relations.  
$$3\alph=3\bet=0,\qquad 3\mu=0,\qquad 3\zeta=0,$$
$$\alph\delt=-\alph\bet,\qquad \delt^2=3\deltwo\bet,
\qquad \alph\deltwo=0,\qquad\delt\deltwo=9\delthr,$$
$$\alph\zeta=0,\qquad\delt\zeta=0,\qquad\alph^2\bet=-\delt\bet^2,$$
$$\alph\mu=0,\qquad\delt\mu=0,$$
$$\deltwo^3-27\delthr^2+\zeta^2=-\delthr(\delt\bet^2+\bet^3)
+\deltwo^2\bet^2+\deltwo\bet^4-(\delt\bet^5+\bet^6).$$

\proof Filtering the ring given by the above relations by powers of
$\bet$ one obtains the ring of Lemma~8, and so if these relations
hold, then they suffice.  It remains to show that these relations do
hold.  The relations $(\delt+\bet)\alph=0$, $\alph\zeta=0$ and
$\alph\deltwo=0$ follow easily using Frobenius reciprocity, because 
$\delt+\bet$, $\zeta$ and $\deltwo$ are defineable as corestrictions
from $\tm$ while $\alph$ restricts to $\tm$ as zero.  The relations
involving $\mu$ follow easily from the fact that the odd degree
cohomology of $\tg$ maps injectively to that of $\tm$.  The
expressions for $\delt^2$, $\delt\deltwo$ and $\zeta\delt$ may also be
shown to hold using Frobenius reciprocity.  As an example, the
following equations verify the expression for $\delt^2$.  
$$\eqalign{\delt^2&=(\Cor^\tg_\tm(\tau)-\bet)\delt
=\Cor^\tg_\tm(\tau)\delt-\delt\bet\cr
&=\Cor(3\tau^2)-\delt\bet=3\deltwo-\delt\bet}$$

From Lemma~8 we see that $\alph^2\bet=a\delt\bet^2+b\alph\bet^2
+c\bet^3$ for some $a$, $b$ and $c$.  By considering the image of this
equation in $\tm$ it may be shown that $c=0$.  As homomorphisms from
$\tg$ to $\tee$, $\alph$ and $\bet$ have kernel containing $\tee$, so
they are in the image of the inflation from $E=\tg/\tee$.  It follows
from Theorem~6 that $\alph^3\bet=\bet^3\alph$.  It is now possible to
solve for $a$ and $b$ using the following equations to show that
either $\alph^2\bet=-\delt\bet^2$ or $\alph^2\bet=\pm\alph\bet^2$.  
$$\bet^3\alph=\alph^3\bet=\alph(a\delt\bet^2+b\alph\bet^2)
=-a\alph\bet^3+b\alph^2\bet^2=(b^2-a)\alph\bet^3+ab\delt\bet^3$$
It remains to rule out the possibility that
$\alph^2\bet=\pm\alph\bet^2$.  The kernel of $\alph+\bet$ viewed as a
map from $\tg$ to $\tee$ is a group isomorphic to $\tp$, and $\alph$, 
$\bet$, $\alph-\bet$ each restrict to this subgroup as a non-zero
element of order three.  Theorem~6 shows that the subring of
$H^*(\tp)$ generated by elements of $H^2$ of order three is isomorphic
to $\Bbb Z[x,y]/(3x,3y,x^3y-y^3x)$.  In this ring the product of any
three non-zero elements of degree two (the degree of $x$ and $y$) is
non-zero.  It follows that $\alph^2\bet-\alph\bet^2$ is non-zero in
$H^*(\tg)$ because its restriction to $\ker(\alph+\bet)$ is non-zero.
Similarly, $\alph^2\bet+\alph\bet^2$ is non-zero because its
restriction to $\ker(\alph-\bet)$ is non-zero.  

There remains now only the expression given for
$\deltwo^3-27\delthr^2+\zeta^2$.  Lemma~8 implies that this quantity
is a multiple of $\bet$, and the relations we have already obtained
show that it is annihilated by $\alph$.  Using Lemma~8, a basis for
$\bet H^{10}$ may be found (this group is isomorphic to $(C_3)^{11}$).
Using the relations already known it may be shown that the kernel of
multiplication by $\alph$ on $\bet H^{10}$ has basis
$\delthr\deltwo\bet$, $\zeta\deltwo\bet$, $\deltwo^2\bet^2$,
$\delthr(\delt\bet^2+\bet^3)$, $\zeta\bet^3$, $\deltwo\bet^4$,
$\delt\bet^5+\bet^6$, and using Proposition~7 it may be shown that
this group maps injectively to $H^{12}(\tm)$.  The relation claimed
now holds, because when multiplied by $\alph$ it gives the valid
relation $0=0$, and its image in $H^{12}(\tm)$ is also a valid
relation.  \qed

\proclaim Corollary 11.  The groups of order 81 are distinguished by
their integral cohomology groups.  

\proof The hardest part of the proof is to distinguish $G(4,1)$,
$G(4,-1)$ and $\gqq$, so we shall only sketch the other cases.  There
are 15 groups of order $3^4$ (see [\burn] for a classification, or
[\thomas] for presentations).  First we use $H^2(G)\cong\Hom(G,\tee)$
to distinguish the five abelian groups, and to split the non-abelian
groups into three classes of sizes three, three and four, with
respective $H^2$ groups $C_3\oplus C_9$, $(C_3)^3$ and $(C_3)^2$.  
The groups $G$ such that $H^2(G)\cong C_3\oplus C_9$ are two split
metacyclic groups, only one of which contains an element of order 27,
and one non-metacyclic.  Three easy spectral sequence arguments show
that the metacyclic with an element of order 27 has $H^3=0$, the other
metacyclic has $H^3\cong C_3$, and the third group has $H^3\cong
C_3\oplus C_3$.  The groups with $H^2\cong (C_3)^3$ are the direct
products of $C_3$ with each of the two non-abelian groups of order 27
and another group isomorphic to the subgroup $P(5,\epsilon)$ of
$G(5,\epsilon)$.  Applying the K\"unneth theorem and Lewis'
description of the cohomology rings of the groups of order $p^3$
[\lew], we see that $H^3\cong (C_3)^4$ for the product of $C_3$ and
the group of order 27 and exponent three, while $H^3\cong (C_3)^2$ for
the other two groups.  These two groups may be distinguished using
$H^4$, because for $P(5,\epsilon)$ this group has exponent nine, while
for the product it has exponent three. 

There remain only the four groups with $H^2\cong(C_3)^2$.  These occur
as normal subgroups of $\tg$ with connected quotient, so we use the
Gysin sequence for $\bee G$ as a $\tee$-bundle over $\bee\tg$ to study
$H^*(G)$.  Using Theorem~10 we may write down a basis for the first
few cohomology groups of $\tg$ as follows.
$$\halign{$#$\hfil&\qquad$#$\hfil\crcr
H^2&\delt,\alph,\bet\cr
H^3&\mu\cr
H^4&\deltwo,\alph^2,\delt\bet,\alph\bet,\bet^2\cr
H^5&\bet\mu\cr
H^6&\delthr,\alph^3,\zeta,\deltwo\bet,\delt\bet^2,\alph\bet^2,\bet^3\cr}$$
Let $\xi$ be one of the elements $\delt-\bet$, $\delt+\bet$,
$\delt+\bet+\alph$, or $\delt$.  (Recall from Proposition~3 that the
kernels of these four elements viewed as maps from $\tg$ to $\tee$
are in distinct $\Aut(\tg)$ classes of normal subgroup, and that they
are isomorphic to $G(4,1)$, $G(4,-1)$, $\gqq$ and the wreath product
of $C_3$ with $C_3$ respectively.)  In each case multiplication by
$\xi$ is injective from $H^2$ to $H^4$, except the case
$\xi=\delt+\bet$, when the kernel has order three.  Thus each of the
four groups has $H^3$ of order three, except $G(4,-1)$ which has $H^3$
of order nine.  In the remaining three cases, $H^4/\xi H^2$ has order
27.  For each of these except the case $\xi=\delt$ multiplication by
$\xi$ is an isomorphism from $H^3$ to $H^5$.  Thus the wreath product
of $C_3$ with $C_3$ has $H^4$ of order 81, while $G(4,1)$ and $\gqq$
have $H^4$ of order 27.  Now consider the remaining two cases.
Multiplication by $\delt+\bet+\alph$ is injective from $H^4$ to $H^6$,
whereas multiplication by $\delt-\bet$ has kernel of order three
(generated by $\delt\bet+\alph^2$).  Hence $H^5(G(4,1))\cong C_3$
whereas $H^5(\gqq)=0$.  \qed

\remarks Lluis and C\'ardenas have found the additive structure of the
cohomology of $C_3$ wreath $C_3$ [\lc], and many authors have studied
the cohomology of the metacyclic groups.  A similar result to
Corollary~11 holds for the groups of order 16, but N.~Yagita has
exhibited a pair of groups of order $p^4$ for all $p\geq 5$ having
isomorphic integral cohomology groups [\yag].  It is not known whether
these groups have isomorphic cohomology rings.  

It will be seen shortly that for $n\geq 5$ the integral cohomology of
$\gne$ is generated by the image of $H^*(\tg)$ and one other element
in degree 5.  The following proposition describes such an element.

\proclaim Lemma 12.  There is an element $\nu$ of $H^5(\gne)$ which
restricts to $M(n,\epsilon)$ as $\gamma\mu$.  The element $\nu$ is not
in the image of $H^5(\tg)$, and may be chosen to be a corestriction
from $H^5(P(n,\epsilon))$.  

\proof Proposition~7 describes the map from $H^*(\tg)$ to $H^*(\tm)$,
and Proposition~4 describes the map from $H^*(\tm)$ to $H^*(M)$.
Together these show that the image of $H^5(\tg)$ in $H^5(M)$ is
trivial, and hence that if $\nu$ exists it cannot be in the image of
$H^5(\tg)$.  To show that $\nu$ exists, we use the double coset
formula $\Res^G_M\Cor^G_P=\Cor^M_N\Res^P_N$.  It is easy to check that
$H^*(N)$ is generated by $\tau$, $\gamma$ and $\mu'$ of degrees 2, 2
and 3 respectively, where $\tau$ and $\gamma$ are the images of the
elements of $H^*(M)$ of the same name, while $\mu'$ satisifes the
equation $\Cor_N^M(\mu')=\mu$ (note that $\Res^M_N(\mu)=0$).  It
follows by Frobenius reciprocity that
$\Cor^M_N(\gamma\mu')=\gamma\mu$, and so it remains to prove that
$\gamma\mu'$ is in the image of the restricion from $P$ to $N$.  The
element $\gamma$ is the image of the element of $H^2(P)$ having the
same name.  To show that $H^3(P)$ maps onto $H^3(N)$, consider the
map of spectral sequences induced by the following commutative diagram.  
$$\matrix{{\bf T}&\mapright{}&BN&\mapright{}&B\tn\cr
\mapdown{{\rm Id}}&&\mapdown{}&&\mapdown{}\cr
{\bf T}&\mapright{}&BP&\mapright{}&B\tp}
$$
For each spectral sequence, $E_3^{3,0}=0$, and the induced map is
surjective on $E_3^{2,1}$, which implies that $H^3(P)$ maps onto
$H^3(N)$.  \qed

\proclaim Theorem 13.  For $n\geq 5$, the integral cohomology ring of
$\gne$ is generated by elements $\alph$, $\delt$, $\mu$, $\deltwo$,
$\nu$, $\delthr$, $\zeta$ of degrees 2, 2, 3, 4, 5, 6 and 6
respectively.  The element $\nu$ is as described in Lemma~12, and the
other generators are the restrictions of the elements of $H^*(\tg)$
(see Theorem~10) having the same names.  They are subject only to the
following relations.  
$$3\alph=0,\qquad 3\mu=0,\qquad 3\nu=0,\qquad 3\zeta=0,$$
$$3^{n-3}\delt=0,\qquad 3^{n-2}\deltwo=0,\qquad 3^{n-1}\delthr=0,$$
$$3\deltwo=\delt^2(1+\epsilon 3^{n-4}),\qquad\delt\alph=0,\qquad
\alph\mu=0,$$
$$9\delthr=\delt\deltwo,\qquad\delt\mu=0,\qquad\deltwo\alph=0,$$
$$\zeta\delt=0,\qquad\zeta\alph=0,\qquad\deltwo^3=27\delthr^2-\zeta^2,$$
$$\delt\nu=0,\qquad\deltwo\nu=\zeta\mu,\qquad
\zeta\nu=-\deltwo^2\mu,\qquad\mu\nu=0.$$

\proof From the Gysin sequence for $\bee\gne$ as a $\tee$-bundle over
$\bee\tg$, whose differential is described in Proposition~3, it
follows that $H^m(\gne)$ is expressible as an extension with kernel
$H^m(\tg)/\xi H^{m-2}(\tg)$ and quotient $\ker(\times\xi :
H^{m-1}(\tg)\rightarrow H^{m+1}(\tg))$, where $\xi =
3^{n-4}\delt-\epsilon\bet$.  To find the kernel of multiplication by
$\xi$, it is helpful to note that the only torsion of order larger
than three in $H^*(\tg)$ is generated by
$\delthr^i(\deltwo^2-3\delthr\delt)$,  which is annihilated by
$\delt$, and so multiplication by $3^{n-4}\delt$ is trivial on the
torsion subgroup of $H^*(\tg)$.  Moreover, it is easy to see that the
kernel of multiplication by $3^{n-4}\delt$ is exactly the torsion
subgroup of $H^*(\tg)$, and hence that $\ker(\times\xi)$ is the
intersection of $\ker(\times\bet)$ and the torsion subgroup of
$H^*(\tg)$.  It may now be checked that $\ker(\times\xi)$ is equal to
the ideal of $H^*(\tg)$ generated by $\alph^2+\delt\bet$ (here it
helps to first find the intersection of $\ker(\times\bet)$ and the
torsion subgroup not in $H^*(\tg)$, but in the ring given by the
filtration of Lemma~8).  It now follows that the element $\nu$ of
Lemma~12 together with the generators for $H^*(\tg)$ forms a
generating set for $H^*(G)$.  

Since $\nu$ was defined as the transfer of an element of order three,
there are no additive extension problems to worry about.  The
relations given in the statement between the generators coming from
$H^*(\tg)$ are just those that hold in $H^*(\tg)$, after substituting
$\epsilon 3^{n-4}\delt$ for $\bet$ throughout.  As a module for the
subalgebra of $H^*(G)$ generated by $\delthr$ and $\alph$, the ideal
generated by $\alph^2+\delt\bet$ is isomorphic to the polynomial
module $\Bbb F_3[\delthr,\alph]$, and the product of
$\alph^2+\delt\bet$ with each of $\delt$, $\mu$, $\deltwo$ and $\zeta$
is zero.  It follows that the only new relations we need to introduce
are expressions for $\delt\nu$, $\mu\nu$, $\deltwo\nu$ and $\zeta\nu$
not involving $\nu$.  Since the restriction from $\tg$ to $\tm$ is
injective in odd degrees, it follows that the restriction from $G$ to
$M$ is injective on the image of $H^\od(\tg)$ in $H^\od(G)$,
and so the expressions given for $\delt\nu$, $\deltwo\nu$
and $\zeta\nu$ may be verified by checking their images in $H^*(M)$
(as described by Propostion~7 and Lemma~12).  The relation $\mu\nu$
follows using Frobenius reciprocity, because $\nu$ is a transfer from
$P(n,\epsilon)$, and the product of any element of $H^5(P)$ and any
element of $H^3(P)$ is zero.  \qed

\proclaim Corollary 14.  For $n\geq 5$, the groups $G(n,1)$ and
$G(n,-1)$ have isomorphic integral cohomology rings.  

\proof  The primed elements of $H^*(G(n,1))$ defined by the equations 
$$\alph'=\alph,\qquad\delt'=\delt,\qquad\mu'=\mu,\qquad
\deltwo'=(1+3^{n-4})\deltwo,$$
$$\nu'=\nu\qquad\zeta'=\zeta,\qquad\delthr'=(1+3^{n-4})\delthr,$$
generate $H^*(G(n,1))$, and satisfy the same relations as the original
generating set for $H^*(G(n,-1))$.  To check that
$\deltwo'^3=27\delthr'^2-\zeta'^2$, it helps to note that the
equations 
$$3\deltwo^2=\deltwo\delt^2(1+\epsilon3^{n-4})=9\delthr\delt$$
imply that $\deltwo^2$ has order $3^{n-4}$. \qed

\remarks Of course, Corollary~14 leaves many questions unanswered,
mainly of the form ``Are ${\cal F}(G(5,1))$ and ${\cal F}(G(5,-1))$
isomorphic?'' for various functors ${\cal F}$.  We list some of these
questions below.  

\def\gone{G(5,1)}
\def\gtwo{G(5,-1)}

Other coefficient rings:  It follows from Corollary~14 and the
K\"unneth theorem that the cohomology groups of $\gone$ and $\gtwo$
are isomorphic for any trivial coefficients.  What can be said about
the ring structure of the cohomology of these two groups with
coefficients $\Bbb Z/(3^m)$?  The strongest result in this direction
would be to show that $\gone$ and $\gtwo$ have isomorphic cohomology
spectra.  (The cohomology spectrum, which was introduced by Bockstein [\bo]
and also studied by Palermo [\pal] consists of the integral and
mod-$m$ cohomology rings for all $m$, together with the projection
maps and Bockstein maps between these rings.)  If the cohomology
spectra of $\gone$ and $\gtwo$ are not isomorphic, this might lead to
examples of pairs of groups distinguishable by their integral
cohomology rings, but not by their integral cohomology groups, using
direct products of copies of $G(5,\epsilon)$.  

Cohomology operations:  One may ask if the isomorphism between the
cohomology rings of $\gone$ and $\gtwo$ commutes with the action of
the algebra of integral cohomology operations.  The author does not
know the answer to this question, but can show that there is an
isomorphism between the two rings commuting with the action of the
operation ``projection to mod-3 coefficients, followed by the first
Steenrod reduced power, followed by the Bockstein back to integer
coefficients''.  This operation is the first possibly non-zero
differential in the Atiyah-Hirzebruch spectral sequence for a 3-group.
One may also ask if Massey products are capable of distinguishing the 
cohomology rings of $\gone$ and $\gtwo$.  Again I do not know the
answer.  

$K$-theory:  The ring $K^0(B\gne)$ is described by Atiyah's
theorem [\ati] together with the presentation for the representation
ring of $\gne$ given in Proposition~2.  I have been unable to decide
whether or not the representation rings of $\gone $ and $\gtwo$ are
isomorphic however.  For any prime $p$, the two non-abelian groups of
order $p^3$ have isomorphic representation rings, and hence isomorphic
$K$-theory rings, although the $\Lambda$-ring structures on their
representation rings are quite different.  This is reflected in the
fact that their cohomology rings are quite different [\lew].  Even including
the $\Lambda$-ring structure, the representation rings of $\gone$ and
$\gtwo$ look very similar.  

Other primes:  As remarked earlier, there is for each odd prime a
family of pairs of groups similar to the groups $\gne$, which may
provide examples for other primes of $p$-groups having isomorphic
integral cohomology rings.  Yagita has been able to show (using these
groups as examples) that there exist non-isomorphic groups of order
$p^4$ for $p\geq 5$ having isomorphic integral cohomology groups [\yag].
Yagita expects to be able to resolve the question of whether these
groups have isomorphic cohomology rings, which may even appear in the
final version of [\yag].  In contrast, work of Rusin [\rus] implies
that even among the groups of order 32 there are no two having
isomorphic integral cohomology.

\par\noindent
{\bf Acknowledgements.}\quad 
I thank my research supervisor, Dr.~C.B.~Thomas, who 
suggested the study of the groups $\gne$ to me and had conjectured
the result of Corollary~9, Dr.~P.H.~Kropholler, who suggested the method
of embedding a finite group in a compact Lie group, and Prof. N.~Yagita 
for helpful correspondence concerning the groups 
$\gne$ and their analogues for 
other primes.  Much of this work was done at Queen Mary and Westfield College, 
where I was supported by SERC post-doctoral fellowship number B90 RFH 8960.

\beginsection References. 

\frenchspacing
\def\book#1/#2/#3/#4/#5/{\item{#1} #2, {\it #3,} #4, {\oldstyle #5}.
\par\smallskip}
\def\paper#1/#2/#3/#4/#5/(#6) #7--#8/{\item{#1} #2, #3, {\it #4,} {\bf #5}
({\oldstyle#6}) {\oldstyle #7}--{\oldstyle#8}.\par\smallskip}
\def\prepaper#1/#2/#3/#4/#5/#6/{\item{#1} #2, #3, {\it #4} {\bf #5} {#6}.
\par\smallskip}

\paper \ati/M. F. Atiyah/Characters and cohomology of finite groups/Publ.
Math. IHES/9/(1961) 23--64/

\paper \bo/M. Bockstein/Homological invariants of the topological product
of two spaces/C. R. (Doklady) Acad. Sci. USSR/40/(1943) 339--342/

\book \burn/W.~Burnside/Theory of Finite Groups/C.U.P./1897/

\paper \lc/H. C\'ardenas and E. Lluis/On the integral cohomology of a
Sylow subgroup of the symmetric group/Comm. Algebra/18/(1990) 105--134/

\book \eve/L. Evens/The cohomology of groups/O.U.P./1991/

\paper \huea/J.~Huebschmann/Perturbation theory and free resolutions for
nilpotent groups of class 2/J. of Algebra/126/(1989) 348--99/

\paper \hueb/J.~Huebschmann/Cohomology of nilpotent groups of class 2/J. of
Algebra/126/(1989) 400--50/

\prepaper \lar/D. S. Larson/The integral cohomology rings of split
metacyclic groups/Unpublished report, Univ. of
Minnesota//({\oldstyle1987})/

\paper \ijlone/I. J. Leary/The integral cohomology 
rings of some $p$-groups/Math. Proc. 
Cambridge Phil. Soc./110/(1991) 25--32/

\prepaper \ijltwo/I. J. Leary/$p$-groups are not determined by their 
integral cohomology groups/submitted//({\oldstyle1992})/

\paper \ly/I. J. Leary and N. Yagita/Some examples in 
the integral and Brown-Peterson cohomology 
of $p$-groups/Bull. London Math. Soc./24/(1992) 165--168/

\paper \lew/G.~Lewis/Integral cohomology rings of groups of order $p^3$/
Trans. Amer. Math. Soc./132/(1968) 501--29/

\paper \pal/F. P. Palermo/The cohomology ring of product complexes/Trans.
Amer. Math. Soc./86/(1957) 174--196/

\paper \rog/K. Roggenkamp and L. Scott/Isomorphisms of $p$-adic group 
rings/Annals of Math./126/(1987) 593--647/

\paper \rus/D. Rusin/The cohomology of groups of order 32/Math. 
Comp./53/(1989) 359--385/

\prepaper \san/R. Sandling/Letter to the author/Jan. 1991///

\book \thomas/C. B. Thomas/Characteristic classes and the cohomology 
of finite groups/Cambridge University Press/1986/

\paper \wei/A. Weiss/Rigidity of $p$-adic $p$-torsion/Annals of 
Math./127/(1988) 317--332/

\prepaper \yag/N. Yagita/Cohomology for groups of ${\rm rank}_pG=2$ 
and Brown-Peterson cohomology/Preprint, 1991///
\end